\documentclass[journal]{IEEEtran}                   

\usepackage[top = 1.0in, bottom = 0.75in, left = 0.75in, right = 0.75in]{geometry}

\usepackage{url}
\usepackage[pdftex]{graphicx}
\DeclareGraphicsExtensions{.eps}
\usepackage[font=footnotesize]{subfig}
\usepackage{caption}
\usepackage{multicol}
\usepackage{algorithm}
\usepackage{algpseudocode}
\usepackage[flushleft]{threeparttable}
\usepackage{wrapfig}
\usepackage{ntheorem}
\usepackage{psfrag}
\usepackage{amsmath}
\usepackage{cite}
\usepackage{epstopdf}
\usepackage{lipsum}
\usepackage{booktabs,caption, dblfloatfix}

\begin{document}

\title{\LARGE{{ Impact of Dynamic Line Rating on Dispatch\\Decisions and Integration of Variable RES Energy}}}


\author{{Bolun Xu,~\IEEEmembership{Student Member,~IEEE}, Andreas Ulbig,~\IEEEmembership{Member,~IEEE}, G\"{o}ran Andersson,~\IEEEmembership{Fellow,~IEEE}}
    \thanks{\footnotesize B.~Xu is with Electrical Engineering, University of Washington. {\textmd Email: xubolun@uw.edu}}%
	\thanks{\footnotesize A.~Ulbig and G.~Andersson are with the Power Systems Laboratory, Dept. of Electrical Engineering, ETH Zurich, 8092 Zurich, Switzerland. {\textmd Email: \{ulbig, andersson\}@eeh.ee.ethz.ch}}
	\vspace{-8mm}
}

\maketitle

\begin{abstract}
	
Dynamic line rating~(DLR) models the transmission capacity of overhead lines as a function of ambient conditions. It takes advantage of the physical thermal property of overhead line conductors, thus making DLR less conservative compared to the traditional worst-case oriented nominal line rating~(NLR).
Employing DLR brings potential benefits for grid integration of variable Renewable Energy Sources~(RES), such as wind and solar energy. In this paper, we reproduce weather conditions from renewable feed-ins and local temperature records, and calculate DLR in accordance with the RES feed-in and load demand data step. Simulations with high time resolution, using a predictive dispatch optimization and the Power Node modeling framework, of a six-node benchmark power system loosely based on the German power system are performed for the current situation, using actual wind and PV feed-in data. The integration capability of DLR under high RES production shares is inspected through simulations with scaled-up RES profiles and reduced dispatchable generation capacity. The simulation result demonstrates a comparison between DLR and NLR in terms of reductions in RES generation curtailments and load shedding, while discussions on the practicality of adopting DLR in the current power system is given in the end.


\end{abstract}

\begin{IEEEkeywords}
	Renewable energy sources, Power generation dispatch, Transmission lines
\end{IEEEkeywords}

\IEEEpeerreviewmaketitle

\section{Introduction}

Facing the challenge of having to reduce $\mathrm{CO_2}$ emissions due to climate change concerns as well as security of supply issues with fossil fuels, many countries nowadays are committed to increasing the share of renewable energy sources (RES) in their electric
power systems, i.e.~wind and PV units. In Germany, for example, the RES share of electricity
generation has increased from 4.7\% of net load demand in 1998 to more than 20\% in 2012. Overall RES electricity generation in 2012 was dominated by wind, PV and hydro generation with an absolute share of net load demand of 8.3\%, 5.0\% and 3.9\%, respectively. The remainder was made up of biomass, land-fill and biogas generation (ca. 3-4\%)~\cite{RESfigure}.

However, existing transmission capacity limitations in many power systems are increasingly impeding the grid integration of ever larger RES energy shares. While building new transmission lines is costly and often requires lengthy legal procedures due to
regulations and public concerns, a short-term alleviation to the capacity limitation problem is to improve the capabilities, and hence the utilization, of existing transmission grids by adopting measures such as dynamic line rating (DLR) for overhead lines.

DLR models the real-time transmission capacity of overhead lines from the conductor thermal balance~\cite{house1958current,morgan1967rating, davis1977new,foss1983dynamic}, hence the current rating is a function of the ambient condition, the physical characteristic of the conductor, and the maximum allowable conductor temperature that protects the line from sag or damaging. 
Foss et al.~\cite{dlr_1990} characterized DLR calculation methods into two approaches: the conductor temperature approach that relies on conductor real-time current and temperature measurements, and the weather approach that calculates DLR from the ambient condition. Because conductor temperature measurements are normally not available, 
the weather approach is more widely used in dispatching studies, because measurements on conductors are not available and DLR is calculated from weather conditions such as air temperature, wind speed, and wind angle~\cite{douglass1996real}-\nocite{ieee_dlr, cigre,wallnerstrom2014impact}\cite{talpur2014implementation}.

Higher wind speed cools down conductors faster, results in a higher transmission capacity for integrating the increased wind generations. By factoring in the cooling effect of wind, a  10-20\% increase in the minimum line rating can be expected in windy areas~\cite{abdelkader2009dynamic}, while a wind speed of 6m/s can at maximum double the line rating compared to the nominal case~\cite{wallnerstrom2014impact}.
A case study showed that applying DLR to the 132kV line between Skegness and Boston enabled the grid integration of 20--50\% more wind generation than by using the more conservative NLR~\cite{dlr_wf}, while another case study applied DLR to a 130~kV regional network and also showed that DLR is significantly profitable in the ampacity upgrading of overhead lines~\cite{talpur2014implementation}.

In this paper, an in-depth investigation of how much DLR can improve the grid integration of RES feed-in in existing power systems
is performed from a grid dispatch prospect with the following focus:

\begin{itemize}
	\item Derive algorithm for DLR model, including line rating and conductor
	surface temperature calculation;
	\item Assess relevant ambient conditions' effect on DLR;
	\item Reconstruction of these ambient conditions;
	\item Establish a benchmark power system model with high renewable
	energy shares and test its improvement on power transmission performance
	once DLR is applied.
	\item Full-year DC OPF dispatch simulations.
\end{itemize}

The paper is structured as follows: Section II introduces the modeling of DLR,
the dispatch model is described in Section III. The design of the benchmark model is explained in Section IV, and the simulation is shown in
Section V, followed by discussion in Section VI and conclusions in Section VII.

\section{Dynamic Line Rating Modeling}

\subsection{Steady-state Heating Balance}

\begin{figure}[!t]
	\centering
	\includegraphics[trim = 09mm 03mm 10mm 03mm, clip, width=.95\columnwidth]{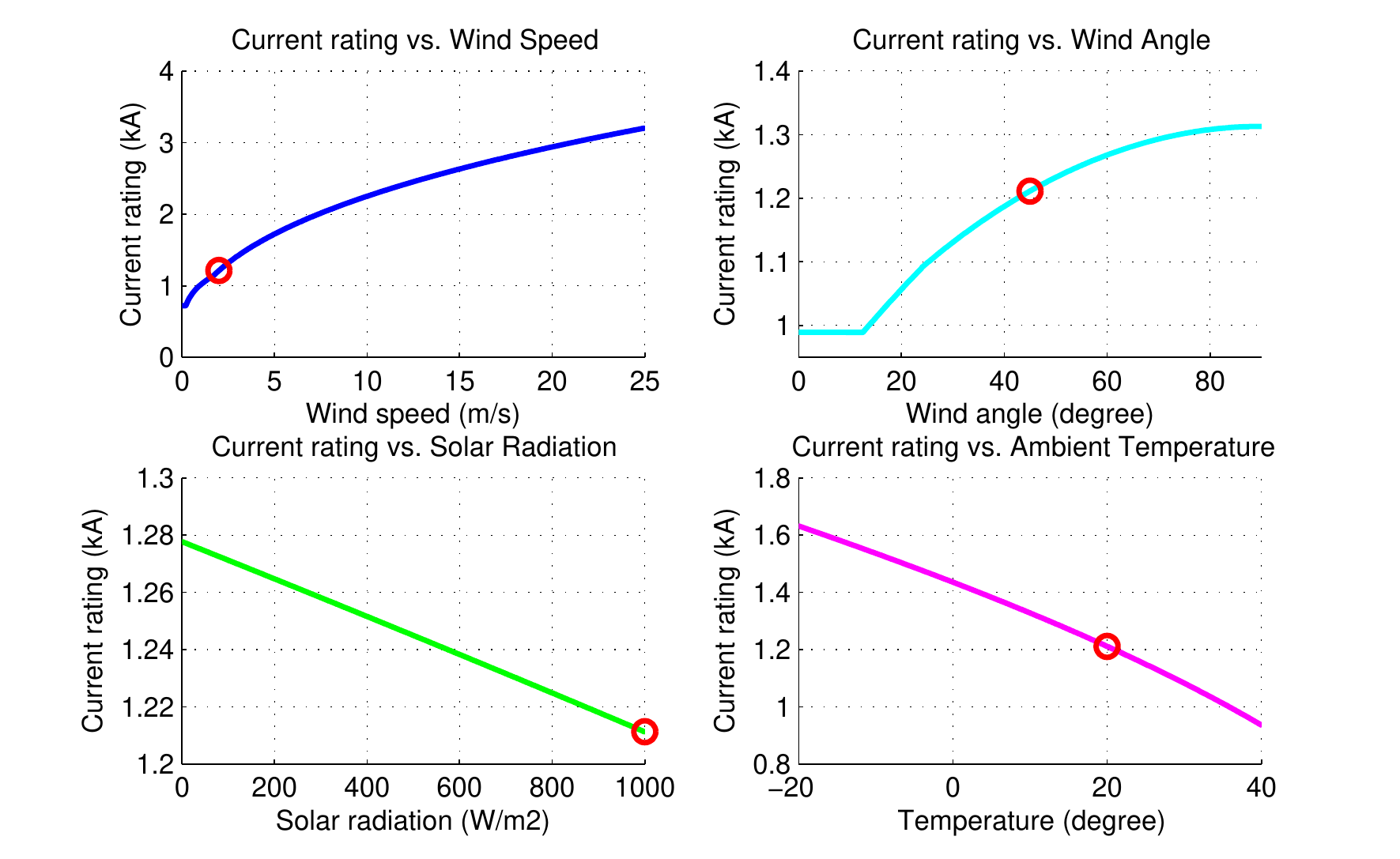}
	\caption{Current rating versus ambient conditions (reference condition marked with red circle).}
	\label{sen:VdST}
\end{figure}

The rating of the transmission line is based on the conductor heat balance in steady-state, defined by a CIGRE standard~\cite{cigre}.
A simplified version of the heat balance equation is

\begin{equation}
	P_{\mathrm{J}} + P_{\mathrm{S}} = P_{\mathrm{c}} + P_{\mathrm{r}} \;,
	\label{bal}
\end{equation}

where $P_{\mathrm{J}}$ is Joule heating, $P_{\mathrm{S}}$ is solar heating, $P_{\mathrm{c}}$ is convective cooling, and $P_{\mathrm{r}}$ is radiative cooling. 

\subsubsection{Joule Heating}
Joule Heating is the resistive heating of conductors. We use the steel cored conductor model in this paper, and calculate $P_{\mathrm{J}}$ as 
\begin{align}
P_J = I^2_{\mathrm{DC}}R_{\mathrm{DC}}\left[1+\alpha(T_{\mathrm{av}}-20^\circ C)\right]\;,
\end{align}
where $R_{\mathrm{DC}}$ is the DC resistance, and $\alpha$ is the resistance temperature coefficient.

\subsubsection{Solar Heating}

Solar heating is due to solar radiations over the conductor surface. $P_{\mathrm{S}}$ depends on the conductor diameter $D$, the conductor surface absorptivity $\alpha_{\mathrm{s}}$, and the global solar radiation $S$, shown as 
\begin{align}
P_{\mathrm{S}} = \alpha_{\mathrm{s}}SD\;.
\end{align}
The value of $\alpha_{\mathrm{s}}$ varies from 0.23 for a bright stranded aluminum conductor, to 0.95 for a weathered conductor in an industrial environment. For most purposes a value of 0.5 may be used.

\subsubsection{Convective Cooling}

The heated surface of the conductor can be cooled down by natural convection (considering no wind) or forced convection (models wind cooling), shown as
\begin{align}
P_{\mathrm{c}} = \pi\lambda_{\mathrm{f}}(T_{\mathrm{s}}-T_{\mathrm{a}})N_{\mathrm{u}}\;,
\end{align}
where $\lambda_{\mathrm{f}}$ is the thermal conductivity of air, $T_{\mathrm{s}}$ and $T_{\mathrm{a}}$ are conductor surface temperature and ambient temperature, respectively. $N_{\mathrm{u}}$ is the Nusselt number, in the natural convection case, $N_{\mathrm{u}}$ is calculated based on conductor surface roughness, and in the forced convection case $N_{\mathrm{u}}$ is calculated from wind velocity and attach angle. The value of $N_{\mathrm{u}}$ used to calculate $P_{\mathrm{c}}$ is the larger one of the two convection cases.

\subsubsection{Radiative Cooling}

Radiative cooling is the cooling due to heat radiations
\begin{align}
P_{\mathrm{r}} = \pi D \varepsilon \sigma_{\mathrm{B}}\left[(T_{\mathrm{s}}+273^\circ C)^4 - (T_{\mathrm{a}}+273^\circ C)^4\right]\;, 
\end{align}
where $\varepsilon$ is emissivity (suggested value is 0.5), $\sigma_{\mathrm{B}}$ is Stefan-Boltzmann constant.

\subsection{Calculate Current Rating in Steady-state Conditions}\label{Sec:DLR_calc}

By assuming coherent temperature distribution across the conductor, the resulting DC current rating($I_{\mathrm{DC}}$) of the heating balance in (\ref{bal}) is


\begin{equation}\label{acdc}
	I_{\mathrm{DC}} = \sqrt{\frac{P_{\mathrm{c}} + P_{\mathrm{r}} - P_{\mathrm{S}}}{R_{\mathrm{DC}}[1+\alpha(T_{\mathrm{av}} - 20)]}} \;,
\end{equation}

where $R_{\mathrm{DC}}$ is the DC resistance, $\alpha$ is the temperature coefficient of the resistance and $T_{\mathrm{av}}$ is the average temperature of the conductor. 
The equivalent AC rating ($I_{\mathrm{AC}}$) can be calculated as

\begin{equation}\label{ac}
	I_{\mathrm{AC}} = \frac{I_{\mathrm{DC}}}{\sqrt{1.0123 + 2.319\cdot10^{-5}I_{\mathrm{DC}}}} \;.
\end{equation}

As shown in Fig.~\ref{sen:VdST},
wind speed ($V$) has a much larger effect on the $I_{\mathrm{AC}}$, increased from $700\textrm{A}$ at $V = 0\textrm{m/s}$ to around $3300\textrm{A}$ at $25\textrm{m/s}$, an increase of 371\%. Besides this, the wind attack angle ($\delta$) and the ambient temperature ($T_{\mathrm{a}}$) also have quite an obvious effect on $I_{\mathrm{AC}}$, with an increase of 35\% and an decrease of 41\%, respectively. The global solar radiation ($S$) has a quite small effect on the rating, and in the simulation $I_{\mathrm{AC}}$ only dropped by 5\% of it's initial value for high solar insolations.

\begin{figure}[h!]
	\centering
	\includegraphics[trim = 5mm 10mm 12mm 10mm, clip, width=.9\columnwidth]{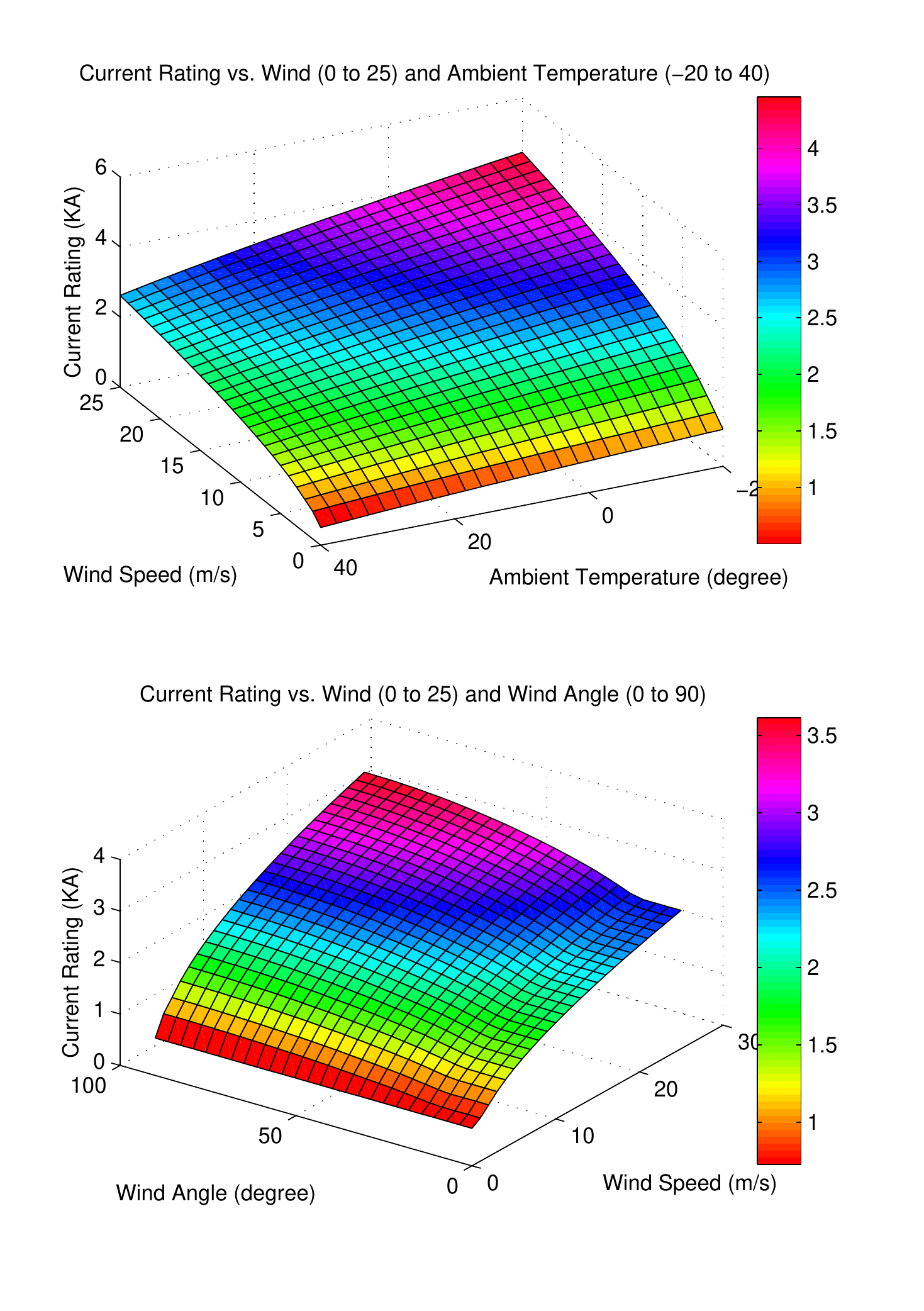}
	\caption{Current rating versus wind attack angle and ambient temperature
	(Default conditions if not specified on axis: $V = 5m/s$, $\delta = 0~to~90^{\circ}$, $S = 1000W/m^2$, $T_a=-20~to~40^{\circ}C$).}
	\label{sen:VdT}
\end{figure} 

\subsection{Calculate Average Conductor Temperature}

The temperature
of the conductor can be calculated from 
$V$, $\delta$, $S$, $T_{\mathrm{a}}$ and current loading $I$
using the following algorithm:

\begin{figure}[!t]
	\centering
	\includegraphics[trim = 40mm 00mm 40mm 00mm, clip, width=.9\columnwidth]{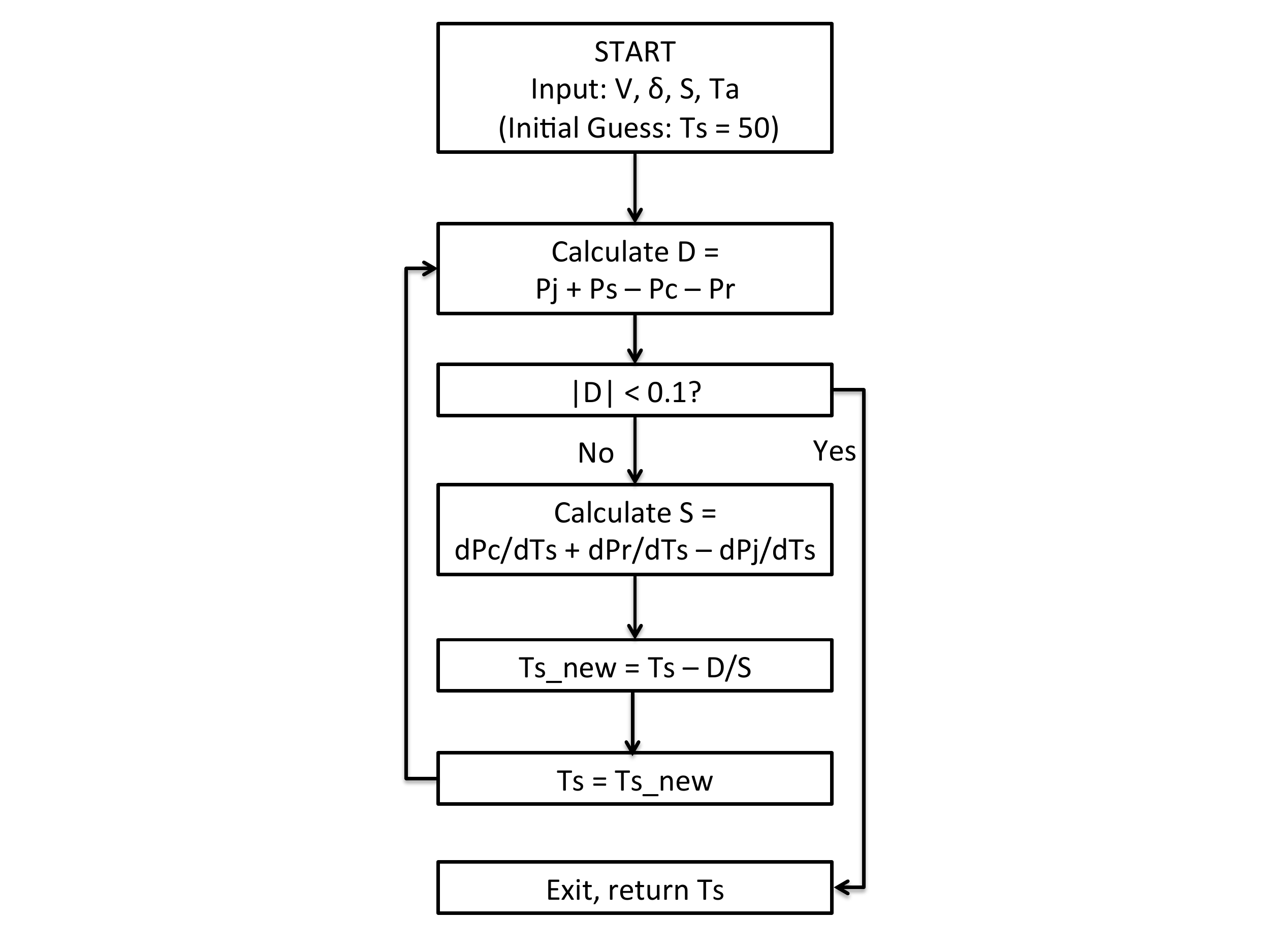}
	\caption{Numerical method for calculating steady-state surface temperature (resulting error $\leq 0.1^\circ C$).}
	\label{sen:ssts}
\end{figure}

Technically, the initial guess of $T_{\mathrm{av}}$ can be any value, however, it makes more
sense to choose a value between the ambient temperature $T_{\mathrm{a}}$ and the maximum
allowable conductor's surface temperature $T_{\mathrm{av, max}}$. From the result of simulation
tests, it is recommended to choose $T_{\mathrm{av}} = 50^{\circ} \textrm{C}$ as a initial guess, with
which the calculation can be finished within three or four iterations. The tolerance is
set to $\delta_{\mathrm{tol}} = 0.1^\circ \textrm{C}$ in this work, which corresponds to an error of 0.014\%
in the calculated DLR value.

\begin{table}[t]
	\begin{center}
		\renewcommand{\arraystretch}{1.0}
		\centering
		\footnotesize
		\caption{Impact of Ambient Parameters on DLR}
		\label{sen:linfit}
		\footnotesize
		\begin{tabular}{lccc}
			\toprule
			\textbf{Situation} & \multicolumn{2}{c}{ \textbf{DLR}} & \textbf{Influence}\\
			&	Absolute	&	Percentage	& \textbf{on DLR} \\
			\toprule
			
			Heating			& $0.14\%$ per 	& $\frac{0.1\%\Delta P_{\mathrm{DLR}}}{1\% \Delta T}$ &			\\
			Load 			& $1^\circ \textrm{C}$  &  &	$7\times$	\\
			Increase		& drop & \scriptsize$T\in[-20,40]^\circ \textrm{C}$&\\
			
			\midrule
			Wind 			& $11.1\%$ per & $\frac{4\%\Delta P_{\mathrm{DLR}}}{1\% \Delta V}$ &\\
			feed-in			& $1\mathrm{m/s}$ 	& &	$300\times$\\
			Increase		& increase & \scriptsize$V\in[0,25] \mathrm{m/s}$&\\
			
			\midrule
			PV 		        & $0.14\%$ per  & $\frac{0.014\%\Delta P_{\mathrm{DLR}}}{1\% \Delta S}$ &	\\
			feed-in			& $100 \mathrm{W/m}^2$  &  & $1\times$	\\
			Increase		& increase & \scriptsize$S\in[0,1000]\mathrm{W/m}^2$ &\\
			
			\bottomrule
		\end{tabular}
	\end{center}
\end{table}

\subsection{Correlation of DLR with Generation \& Load Volatility}

While wind power feed-in is proportional to $V^3$ and
solar power feed-in is proportional
to $S$ \cite{masters},
there exists a well-known negative proportionality of load demand and ambient,
i.e. outdoor, temperature during the winter season, such as in France \cite{euro_heat}.
A linearized model which relates DLR with generation and consumption is
presented in Table~\ref{sen:linfit}, which is by fitting a first-order
curve to the steady-state rating analysis in the previous section (wind
angle is assumed to be $45^\circ$ and is not included in this table).
From the percentage result we can see that wind has an influence that is
40x stronger than the influence of ambient temperature on DLR ($\frac{4.0\%/\%}{0.1\%/\%} = 40$),
and close to 300x stronger than the influence of solar insolation
($\frac{4.0\%/\%}{0.014\%/\%} \approx 286$).

\section{Description of Simulation \& Analysis}

With the calculation methods for DLR established, it is now possible to add DLR to a power dispatch model and examen the improvement of the dispatch result. In this work, a previously established economic power dispatch model is used \cite{power_node, jonas, fortenbacher}, in which power flows are calculated using the DC approximation method, and dispatchable and nondispatchable generators, controllable and non-controllable loads as well as different types of storage units are modeled and processed.

\subsection{Economic Dispatch Model}

The economic dispatch optimization problem in time step k with
objective function $J(k)$ can be formulated as follows:
\begin{align}
	\min J(k) = & \sum_{l=k+N-1}^{l=k} (x(l)-x_{ref})^T\cdot Q_x \cdot (x(l)-x_{ref})\nonumber\\
	 & R_u \cdot u(l) +\delta u(l)^T \cdot \delta R_u \cdot \delta u(l) \nonumber\\
	s.t. \quad 
	(a) \quad & x(l+1) = A \cdot x(l) + B \cdot u(l) \nonumber\\
	(b) \quad & 0\leq x^{min} \leq x(l) \leq x^{max} \leq 1 \nonumber\\
	(c) \quad & u^{min} \leq u(l) \leq u^{max}\nonumber\\
	(d) \quad & \delta u^{min} \leq \delta u(l) \leq \delta u^{max}\nonumber\\
	(e) \quad & P_{line}^{min}(l) \leq P_{line}(l) \leq P_{line}^{max}(l)\quad,
\end{align}

where $N$ is the prediction step number, $x$ and $x_{ref}$ are state variables and their reference value,
$u$ is the node variables, and $Q_x$, $R_u$ and $\delta R_u$ are optimization parameters. 
The sampling time is 15 minutes, and the prediction horizon is 64 hours, with perfect prediction, i.e. accurate load, wind and PV forecasts, assumed.
Equation (a) representing linear Power Node equations and (b-e) are system constraints.
Notably,  topology constraints $P_{line}^{min}$ and $P_{line}^{max}$ are assumed to be constants values in the case of NLR, whereas in the case of DLR the line limits become time variant, i.e. $P_{line}^{min}(k)$ and $P_{line}^{max}(k)$.

\subsection{Simulation Framework}

\begin{figure}
	\centering
	\includegraphics[trim = 80mm 00mm 20mm 00mm, clip, width=.9\columnwidth]{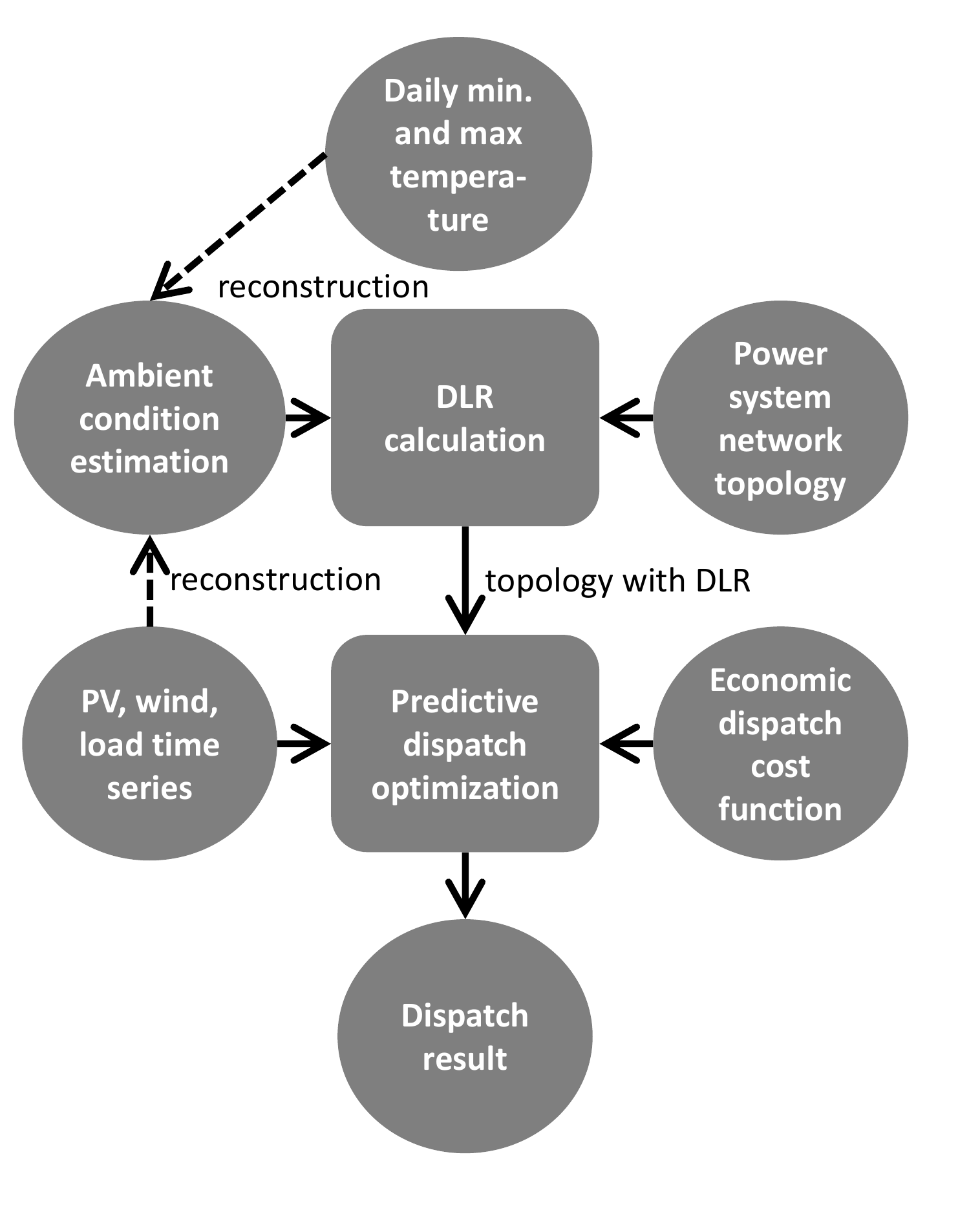}
	\caption{Economic power dispatch simulation framework with DLR}
	\label{opt:DLR}
\end{figure}

The economic dispatch simulation set-up is shown in Fig.~\ref{opt:DLR}. In this project, the ambient conditions are estimated from wind and PV feed-in series and daily temperature records, detailed estimation procedures are later illustrated in Section~\ref{Sec:rec_amb}. The DLR of each line in the system topology are calculated from the ambient condition in its region and the conductor parameters. Dispatch decisions are made based on system line ratings (can be DLR or NLR),  renewable feed-in and load series, and the economic dispatch cost function as described in the previous section. 

\section{Benchmark System \& Simulation}

Germany is chosen as prototype reference for the benchmark model used for the dispatching simulation.
In Germany, wind is
stronger in the north while solar insulation is stronger in the south, most
wind turbines are installed in the northern part, and the majority of
the PV units are installed in the south \cite{energymap}. So if Germany
is to improve its share of RES power, it will be facing the problem to
transmit the wind and solar power nationwide, especially the transmission capacity
bottleneck in the north-south direction~\cite{denaII},
which makes it a fine
model to test the performance of DLR.

\begin{figure}[t]%
\centering
\subfloat[TSO zones]{
	\includegraphics[trim = 0mm 0mm 0mm 0mm, clip, width = 0.45\columnwidth]{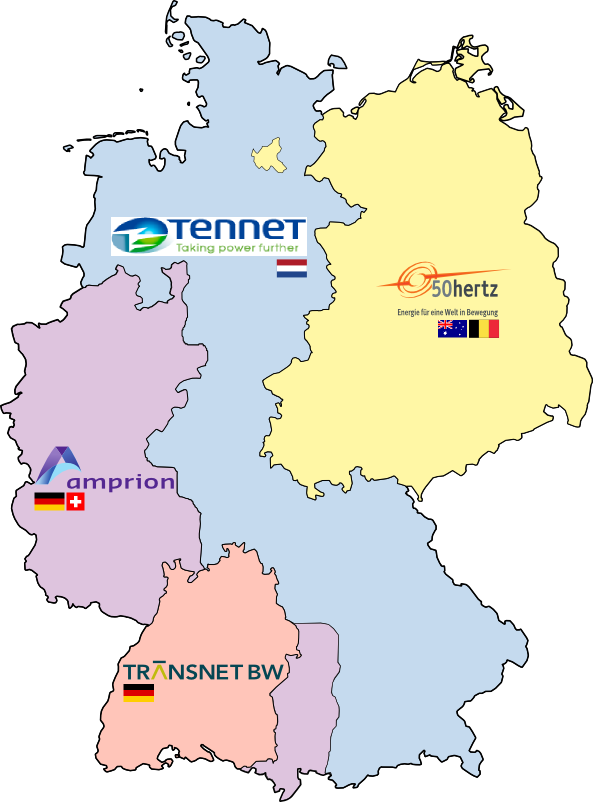}
	\label{fig:benchmarka}%
}%
\subfloat[Benchmark zones]{
	\includegraphics[trim = 0mm 0mm 0mm 0mm, clip, width = 0.45\columnwidth]{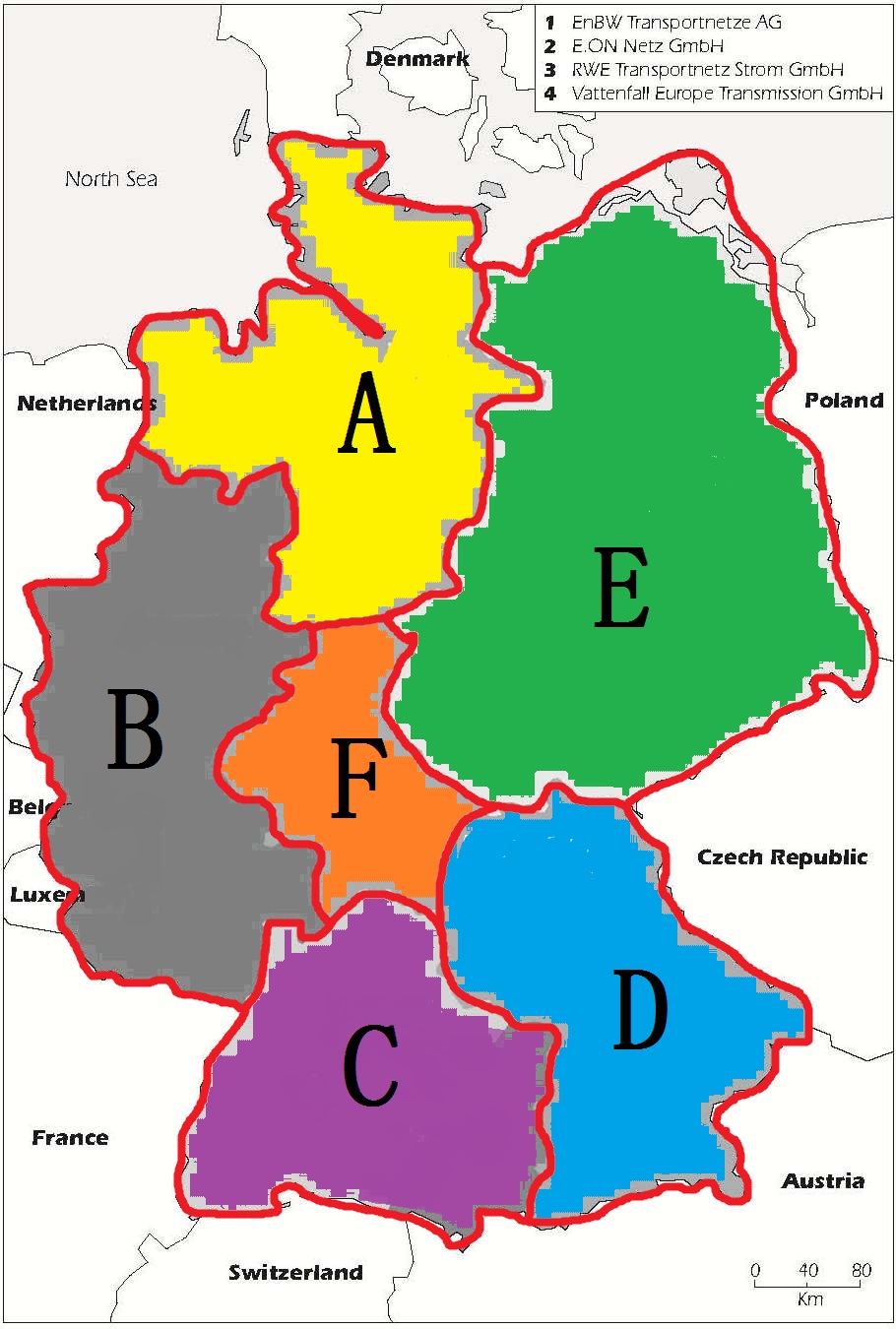}
	\label{fig:benchmarkb}%
}
\caption{The 6-node benchmark model based on TSO zones in Germany.}%
\label{fig:benchmark}%
\end{figure}

\begin{figure}[t]
	\centering
	 \includegraphics[trim = 0mm 10mm 00mm 5mm, clip, width=.9\columnwidth]{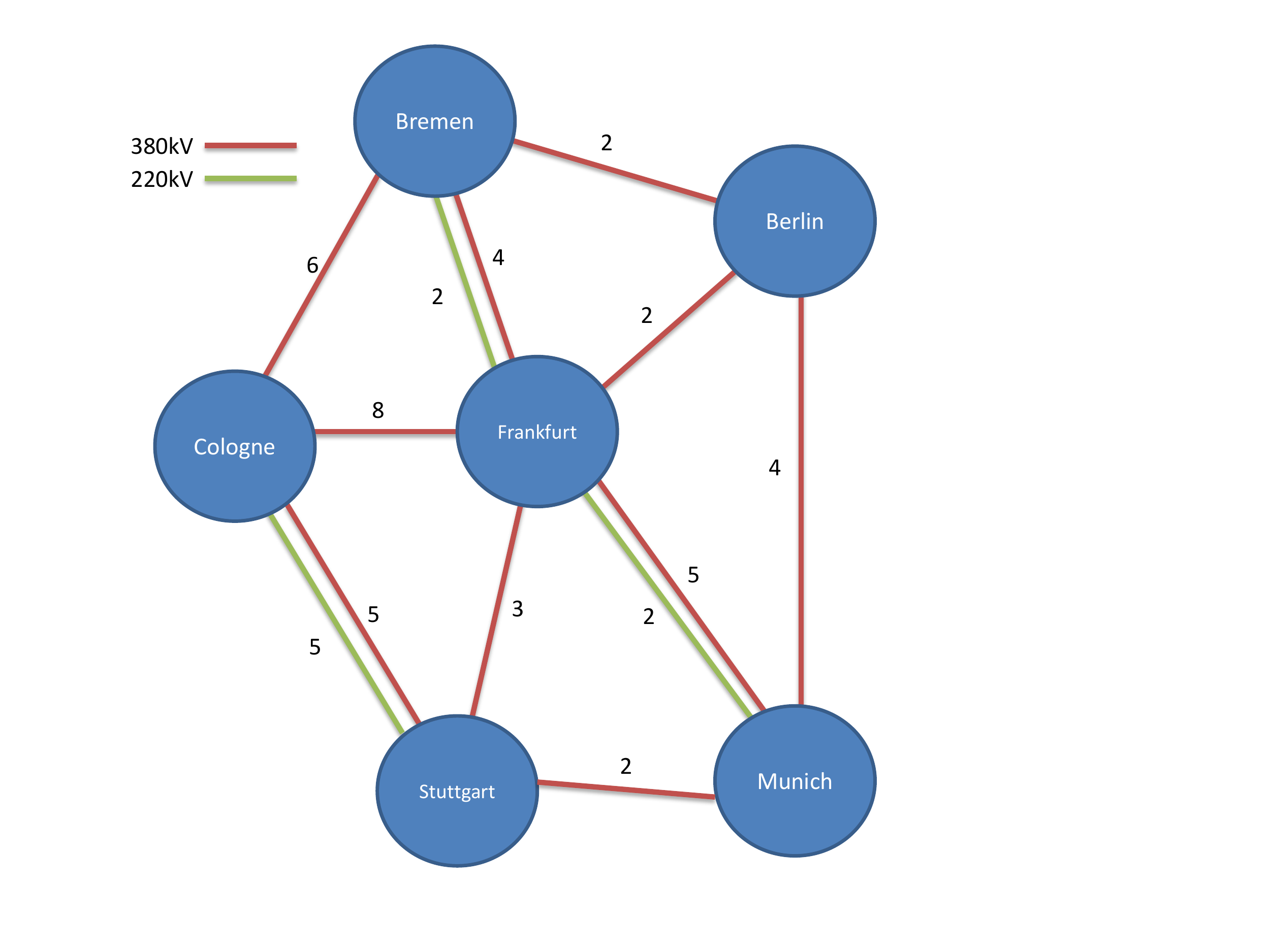}
	 \caption{Simplified 6-node grid topology (red lines indicate number of 330kV transmission lines, green lines 220kV).}
	 \label{fig:bench_topo}
\end{figure}

\subsection{Topology Design of the Benchmark Model}

The four transmission system operators (TSO) in Germany (Fig. \ref{fig:benchmarka}) roughly split Germany's transmission network into four zones \cite{tso_germany}\cite{tso_wiki}. While having other TSO zones unmodified, the Tennet TSO zone is split into three smaller zones, in the order of north, middle and south. The border for these three zones are made where high-voltage overhead lines cross and no distribution network exists.

Based on the assumption that the power feed-in and consumption over one
zone can be simplified into a single node (city), the weather condition
over one zone is the same and can be represented by the weather data of
this city, the six zones can be simplified into
six nodes, represented by six cities: Bremen (A), Cologne (B), Stuttgart
(C), Munich (D), Berlin (E) and Frankfurt (F). The transmission capacity
of two zones is determined by counting the number of 220kV and
380kV overhead
lines connecting them, using the ENTSO-E network
map~\cite{entso_e_map}.
Fig.~\ref{fig:benchmarkb} shows the designed benchmark model.
Fig.~\ref{fig:bench_topo} shows the simplified 6-node grid topology.

\begin{table}[t]
	\begin{center}
	\renewcommand{\arraystretch}{1.0}
	\centering
	\footnotesize
	\caption{Federal states contained in benchmark zones}
	\label{tab:ben:state}
	\begin{tabular}{lc}
	\toprule
	Benchmark Zone		& 	Corresponding Federal States in Germany				\\
	\midrule
	A \scriptsize{(BREMEN)}			& Bremen, Hamburg, Niedersachsen, 						\\
									& Schleswig-Holstein									\\
	\midrule
	B \scriptsize{(COLOGNE)}		& Nordrhein-Westfalen, Rheinland-Pfalz, Saarland		\\
	\midrule
	C \scriptsize{(STUTTGART)}		& Baden-W\"{u}rttemburg, Bayern (south-east part)								\\
	\midrule
	D \scriptsize{(MUNICH)}			& Bayern (all except south-east part)												\\
	\midrule
	E \scriptsize{(BERLIN)}			& Berlin, Brandenburg, Mecklenburg-Vorpommern,  		\\
									& Hamburg, Sachsen, Sachen-Anhalt, Th\"{u}ringen					\\	
	\midrule
	F \scriptsize{(FRANKFURT)}		& Hessen												\\
	\bottomrule
	\end{tabular}
	\end{center}
\end{table}

\subsection{Ambient Data Reconstruction}\label{Sec:rec_amb}

Calculating DLR values in accordance with the 15-minutes simulation time step requires wind velocity, solar radiation, and ambient temperature data in the same time step. Obtaining such historical measurements in the six designed zones in Germany can be a tremendous work and is beyond the effort of this study. As an alternative approximation, wind velocity and solar radiation data are reconstructed from wind and PV feed-ins, while ambient temperature data are generated from daily maximum and minimum temperature with a sinusoidal approximation method. 

\subsubsection{Wind Velocity}

The electricity power produced by a wind turbine ($P_{\mathrm{W}}$) can be described as 
\begin{equation}\label{rec:wind}
	P_{\mathrm{W}} = \frac{1}{2}C_{\mathrm{p}}A\rho V^3\;,
\end{equation}
where $A$ is the rotor area of a wind turbine, $\rho$ is the air density, and $C_{\mathrm{p}}$ is the coefficient of performance \cite{masters2013renewable}. By assuming the wind velocity is always between the rated wind speed ($V_{\mathrm{rated}}$) and the cut-in ($V_{\mathrm{cut}}$) wind speed, and assume $C_{\mathrm{p}}$ be a constant for all wind speeds, then the wind power generation is proportional to $V^3$. We further assume that $A$ is the same for all wind turbines in Germany and $\rho$ is the same over all Germany, and the number of wind turbines operating in Germany is fixed during the year 2011. Therefore the highest feed-in power ($P_{\mathrm{W,max}}$) can be mapped to the rated wind speed and the minimum wind power ($P_{\mathrm{W}}=0$) to the cut-in wind speed, define $C_{\mathrm{W}} =  \frac{1}{2}C_{\mathrm{p}}A\rho$, then it follows that
\begin{align}
C_{\mathrm{W}} = \frac{P_{\mathrm{W,max}}}{(V_{\mathrm{rated}}-V_{\mathrm{cut}})^3}\;.
\end{align}

According to the wind turbine design from ENERCON Inc.~\cite{enercon}, rated wind speed is set to 15m/s and cut-in wind speed is 1m/s. This is a range that wind speeds will stay in with a probability of 95\%~\cite{ebuchi2005intercomparison}. Thus, the wind speed can be calculated from the corresponding wind feed-in as
\begin{align}
V = \sqrt[3]{\frac{P_{\mathrm{W}}}{C_{\mathrm{W}}}} + V_{\mathrm{cut}}\;.
\end{align}

\subsubsection{Solar Radiation}

The current and voltage generated by a photovoltaic cell can be represented as \cite{masters2013renewable}
\begin{align}
I = I_{\mathrm{SC}} - I_{\mathrm{d}}\;,\quad V = V_{\mathrm{d}}-IR_{\mathrm{S}}\;,
\end{align}
where $I_{\mathrm{SC}}$ is the short circuit current and is proportional to the solar radiation, $R_{\mathrm{S}}$ is the series resistance, $I_{\mathrm{d}}$ and $V_{\mathrm{d}}$ are the PV diode's current and voltage. By assuming $R_{\mathrm{S}}$, $I_{\mathrm{d}}$, and $V_{\mathrm{d}}$ are constants, the power delivered by a PV cell is  proportional to the solar radiation. 
We assume the same integrated PV capacity in Germany throughout a year, and map the PV feed-in power to the solar radiation data using a linear correlation
\begin{align}
S = \frac{S_{\mathrm{mean}}}{P_{\mathrm{PV,mean}}}P_{\mathrm{PV}}\;
\end{align}
where $P_{\mathrm{PV,mean}}$ and $S_{\mathrm{mean}}$ are the average PV feed-in power and the average solar radiation in the region of interest.

\subsubsection{Ambient Temperature}

The temperature variation within a day can be represented with a sinusoidal approximation as \cite{ephrath1996modelling}
\begin{equation}\label{rec:temp_sin}
T_{\mathrm{a}}(t) = T_{\mathrm{min}} + (T_{\mathrm{max}} - T_{\mathrm{min}})\Gamma(t)\;,
\end{equation}
where $T_{\mathrm{a}}$ is the air temperature, $T_{\mathrm{max}}$ and $T_{\mathrm{min}}$ are the maximum and minimum temperature in a day, respectively. $\Gamma(t)$ is the sinusoidal approximation function, ranging from 0 to 1. $\Gamma(t)$ is designed to match its peaks and valleys to the occur time of $T_{\mathrm{min}}$ and $T_{\mathrm{max}}$ in each day, which at $45^\circ$ north latitude are approximately 4 and 18 o'clock in summer seasons, 8 and 14 o'clock in winter seasons, and 6 and 16 o'clock in spring and autumn seasons~\cite{pidwirny2006daily}.

\subsection{RES Power feed-in and Load Demand Data}

The generation and load profile with high time-resolution (15 minutes) are obtained from the
four TSOs in Germany. Profiles from Amprion, EnBW, and 50Hertz are mapped to benchmark zone
B, C, and E respectively. Load profile from  TENNET are divided to benchmark zone A, D, and F according to population proportion, while wind and PV generation profile are divided according to unit installation capacity proportion. With Table \ref{tab:ben:state}, the Pumped-Storage Hydroelectricity (PSH) capacity is also determined using the installed PSH capacity in each federal state \cite{pump}. The dispatchable generation capacity is set to 78 GW in total and is split according to the share of population \cite{energymap}.



\subsection{Daily Maximum and Minimum Temperature}

Building on a traditional dispatch model (NLR dispatch model), the purposed DLR dispatch requires no additional data except the daily maximum and minimum temperature records. In the designed benchmark power system model, the temperature records in the six zone cities are used to represent the weather over the entire zone, and are obtained from the German Federal Ministry of Transport, Building and Urban Development \cite{dwd}.
\begin{figure*}[h]
	\centering
	\vspace{-10mm}
	\includegraphics[trim = 20mm 10mm 20mm 05mm, clip, width = 1.7\columnwidth]{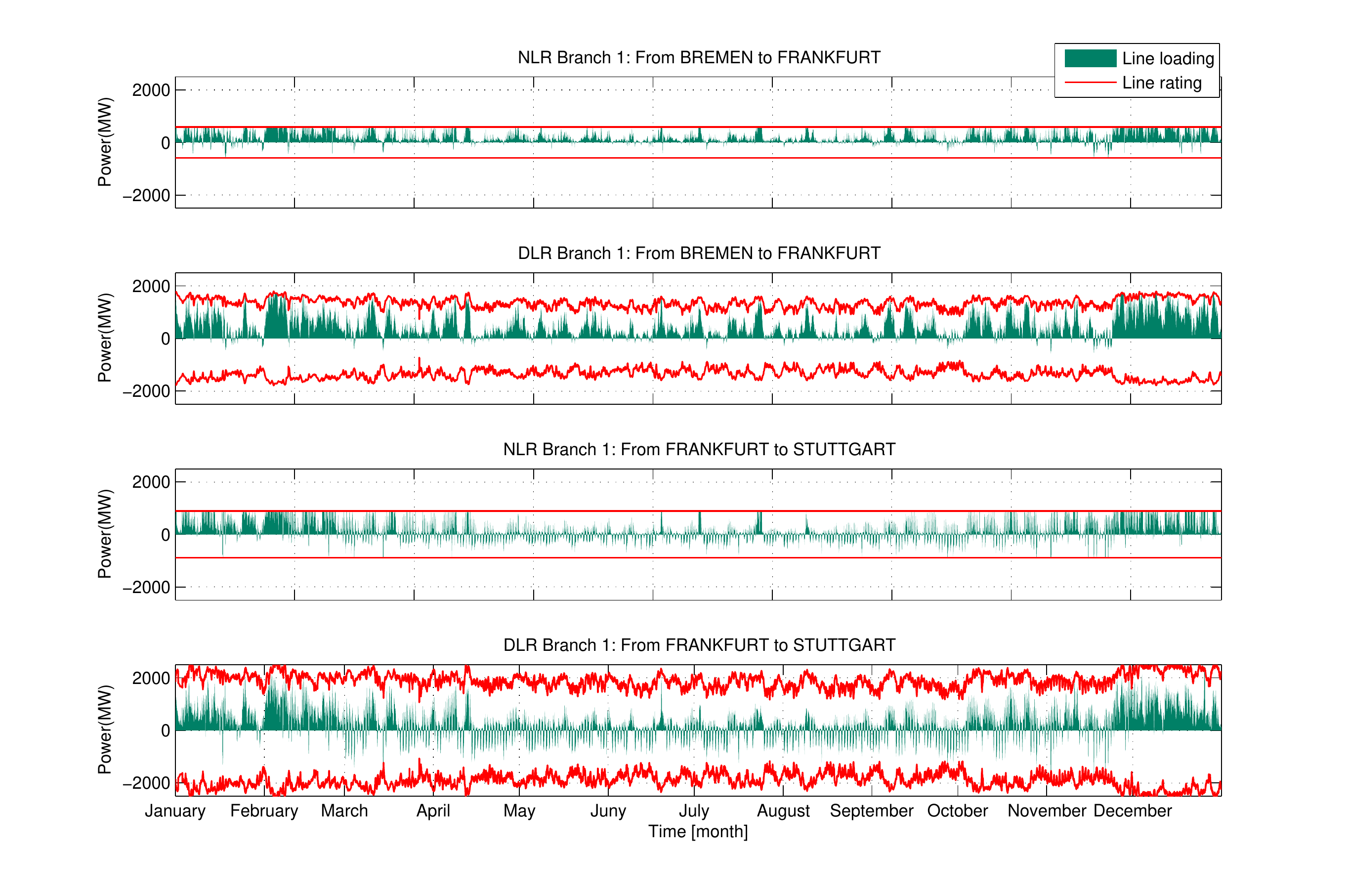}
	\caption{Transmission line loadings of BREMEN-FRANKFURT and FRANKFURT-STUTTGART.}
	\label{fig:NLR_DLR}
\end{figure*}

\begin{figure*}[t]%
\centering
\subfloat[BREMEN zone NLR]{
	\includegraphics[trim = 60mm 103mm 120mm 98mm, clip, width = .95\columnwidth]{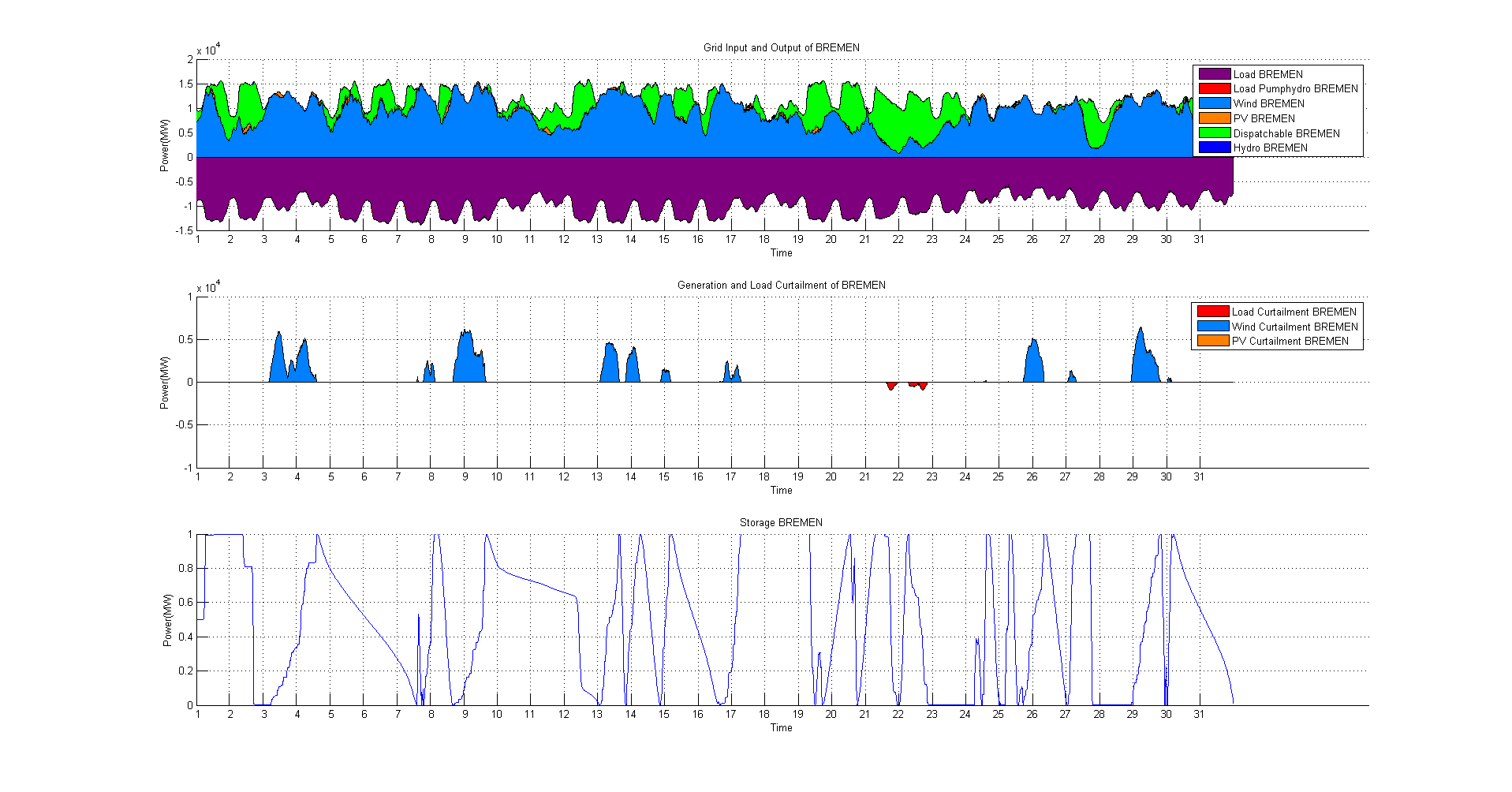}
	\label{fig:bremen_nlr}%
}%
\subfloat[BREMEN zone DLR]{
	\includegraphics[trim = 60mm 103mm 120mm 98mm, clip, width = .95\columnwidth]{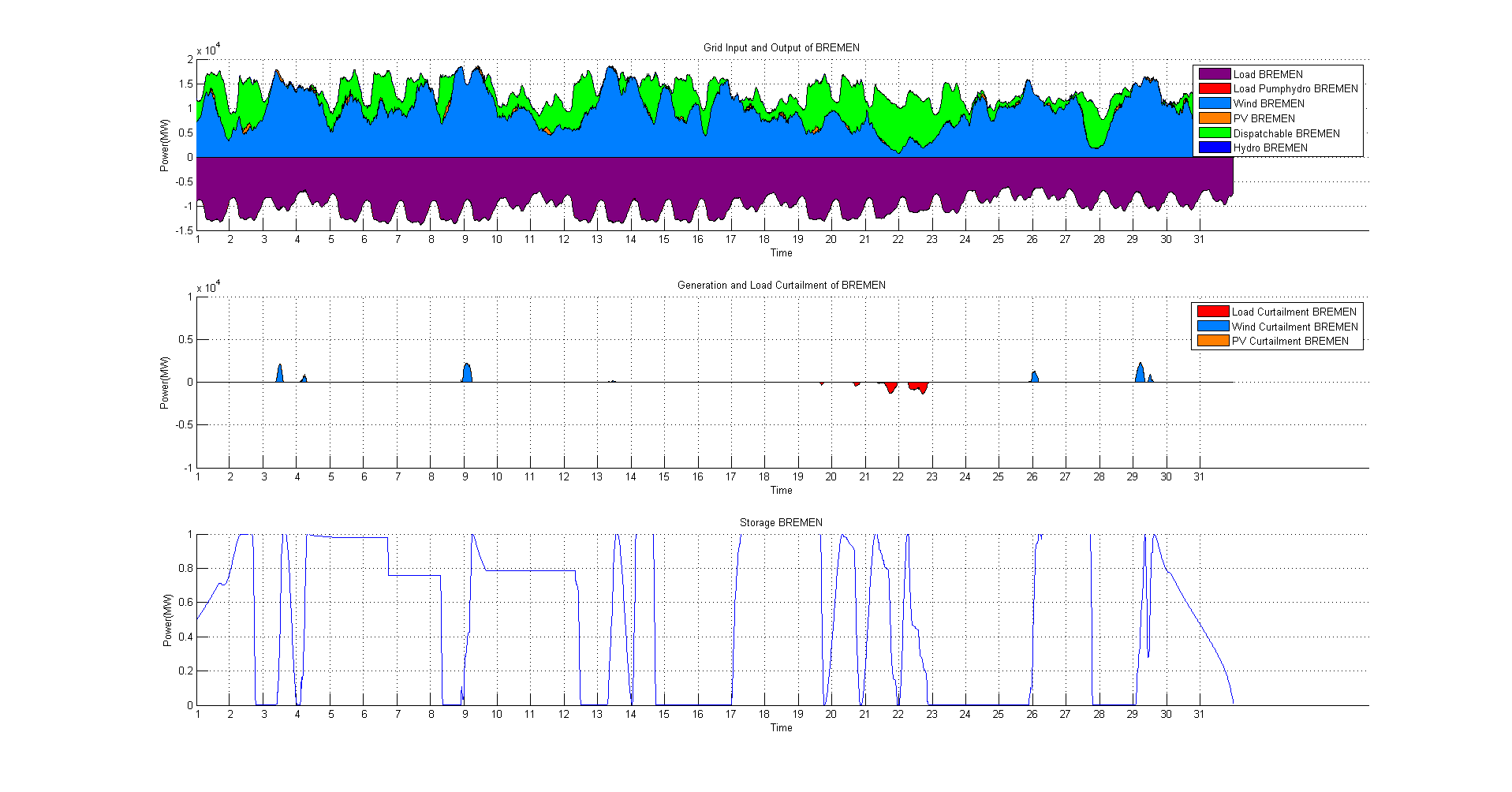}
	\label{fig:bremen_dlr}%
}
\\
\subfloat[STUTTGART zone NLR]{
	\includegraphics[trim = 60mm 103mm 120mm 98mm, clip, width = .95\columnwidth]{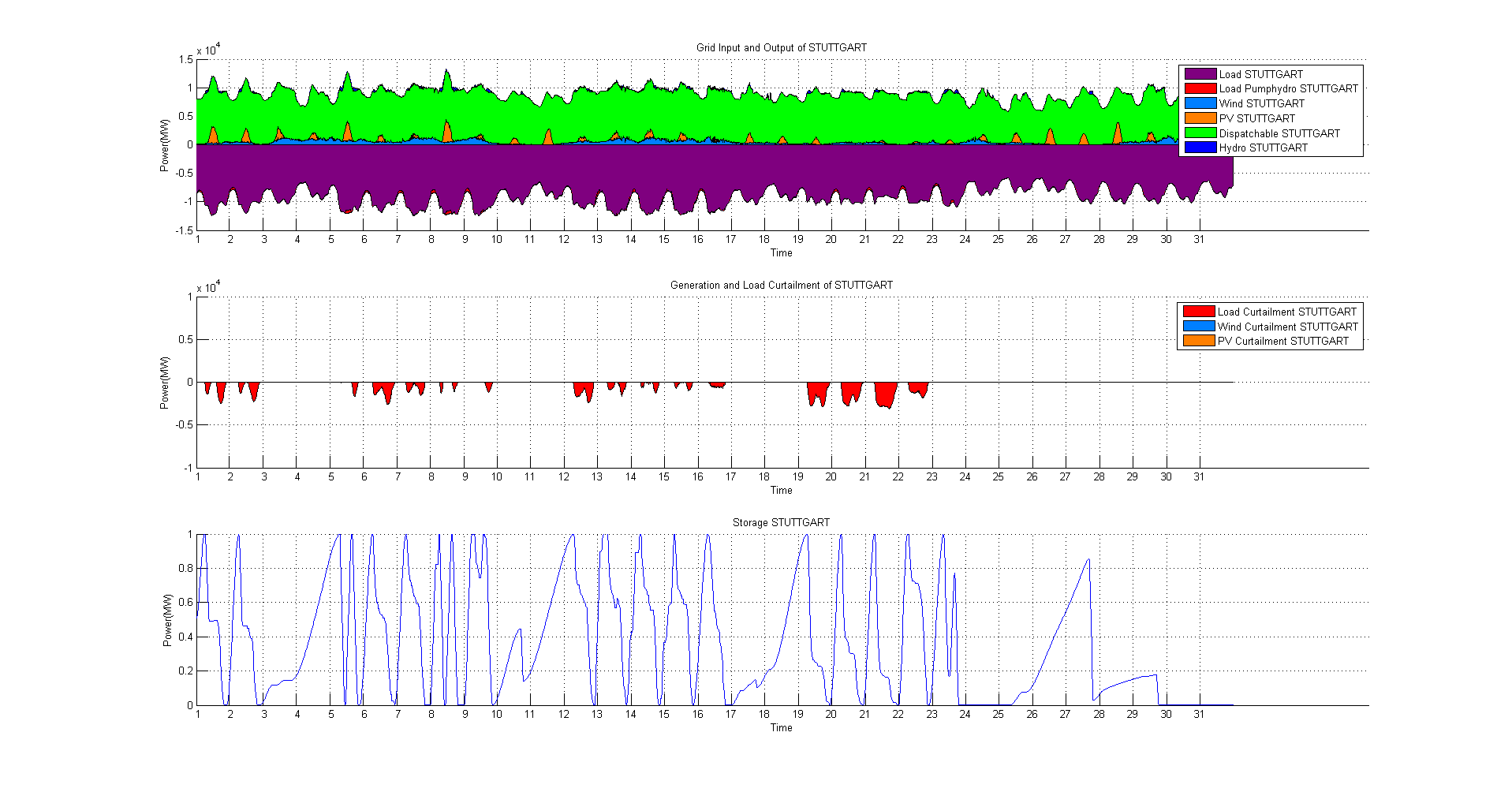}
	\label{fig:stuttgart_nlr}%
}%
\subfloat[STUTTGART zone DLR]{
	\includegraphics[trim = 60mm 103mm 120mm 98mm, clip, width = .95\columnwidth]{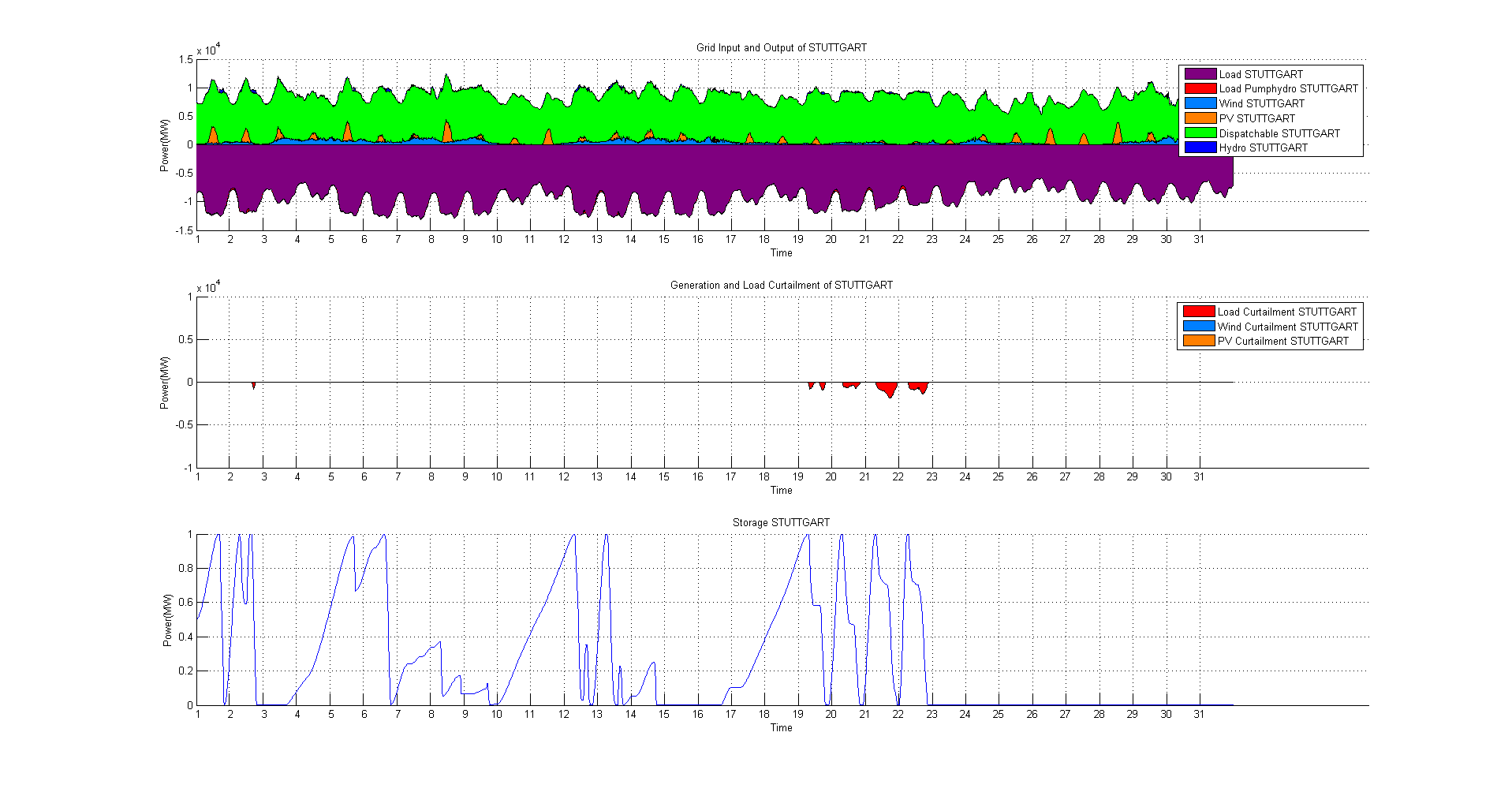}
	\label{fig:stuttgart_dlr}%
}
\caption{Comparison of hourly renewable curtailment (positive areas in blue) and load shedding (negative areas in red) in December (the unit of X axis is day in the month; the unit of Y axis is 10~GW, ranging from 10~GW to -10~GW).}%
\label{fig:bremen_stuttgart}%
\end{figure*}


\subsection{Calculation of Dynamic Line Rating (DLR)}

\begin{table}[t]
	\begin{center}
		\renewcommand{\arraystretch}{1.0}
		\centering
		\footnotesize
		\caption{NLR vs. DLR comparison (2011)}
		\label{tab:rating}
		\begin{footnotesize}
			\begin{tabular}{lcccc}
				\toprule
				
				Zone			& NLR (MVA) 	& \multicolumn{3}{c}{DLR (MVA)}		\\
				Connection		& 	& min	& mean	& max \\
				\midrule
				A $\leftrightarrow$ B				& 1778		& 1686		& 3144		& 4128	\\
				\midrule
				A $\leftrightarrow$ E			& 1529		& 1470		& 2721		& 3572	\\
				\midrule
				A $\leftrightarrow$ F		& 592		& 560		& 1052		& 1382	\\
				\midrule
				B $\leftrightarrow$ C		& 2340		& 2029		& 3884		& 5456	\\
				\midrule
				B $\leftrightarrow$ F		& 2371		& 2189		& 4109		& 5515	\\
				\midrule
				C $\leftrightarrow$ D		& 592		& 454		& 948		& 1397	\\
				\midrule
				C $\leftrightarrow$ F		& 889		& 819		& 1482		& 2041	\\
				\midrule
				D $\leftrightarrow$ E		& 1186		& 1119		& 1990		& 2775	\\
				\midrule
				D $\leftrightarrow$ F		& 1825		& 1692		& 3052		& 4231	\\
				\midrule
				E $\leftrightarrow$ F		& 583		& 555		& 1034		& 1399	\\
				\bottomrule
			\end{tabular}
		\end{footnotesize}
	\end{center}
\end{table}

Because DLR is derived from the thermal balance of conductors, calculating DLR is equivalent to finding the current magnitude that generates the maximum allowable temperature under the weather condition. We define the maximum allowed conductor temperature to be $75^\circ C$, and assume the conductor type of all HV transmission lines is 428-A1/S1A-54/7 'Zebra'. Minor scalings are made to fit the DLR result to the corresponding NLR in the system. We also assume the wind attack angle being fixed at $45^\circ$. The ambient condition estimations are made to the six benchmark zones, while all transmission lines are the interconnection between two zones, thus for a transmission line, its DLR-based transmission line ratings under the ambient condition in each of the two connected zones are calculated, and the lower value is used.

In Table~\ref{tab:rating}, the estimated NLR and DLR of the benchmark model are shown. It can be seen that the minimum value of DLR is slightly lower than the NLR value. NLR is normally determined under $35^\circ$ air temperature, $45^\circ$ wind attack angle, $1000 \mathrm{W/m^2}$, and a wind speed around 1m/s. Thus for extreme events with no wind and high air temperature, the calculated DLR values can be lower than NLR values. 

\section{Results}


\begin{table}[t]
	\begin{center}
		\renewcommand{\arraystretch}{1.0}
		\centering
		\footnotesize
		\caption{Scaled-up Dispatch Simulation Result (full-year 2011, in percentage of the total loads or generations)}
		\label{tab:result}
		\begin{footnotesize}
			\begin{tabular}{lcccccc}
				\toprule
				Benchmark			&  \multicolumn{3}{c}{NLR Curtailments (\%)}	&  \multicolumn{3}{c}{DLR Curtailments (\%)}	\\
				Zone 				& Load & Wind & PV 	& Load & Wind & PV\\
				\midrule
				BREMEN 		& 0.366 & 2.209 & 0 	& 0.429 	& 0.183 & 0\\
				\midrule
				COLOGNE 	& 0.234 & 0 	& 0 	& 0.245 	& 0 	& 0\\
				\midrule
				STUTTGART 	& 0.804 & 0 	& 0 	& 0.457 	& 0 	& 0\\
				\midrule
				MUNICH 		& 0.428 & 0 	& 0 	& 0.529 	& 0 	& 0\\
				\midrule
				BERLIN 		& 0.274 & 0.534 & 0 	& 0.339 	& 0 	& 0\\
				\midrule
				FRANKFURT 	& 1.826 & 0 	& 0 	& 1.012 	& 0 	& 0\\
				\midrule
				TOTAL 		& 0.493 & 1.086 & 0 	& 0.420 	& 0.047 & 0\\
				\bottomrule
			\end{tabular}
		\end{footnotesize}
	\end{center}
\end{table}

Simulations are performed on the proposed benchmark model with actual or scaled renewable feed-in and load demand time-series profile with a sampling time of 15 minutes for the full-year 2011. The rating of the transmission lines can be either the NLR or DLR value, and the wind and PV generation curtailment and the load shedding of the system is then compared. 

\subsection{Renewable Generation and Load Curtailment}

The simulation with original profiles showed no curtailment in both NLR and DLR simulations. In the scaled-up case, profiles are scaled by doubling wind and PV feed-ins and decreasing dispatchable generation capacities to 80\% of the original value, resulting in abundant wind resources in the north and generation shortages in the south, creating an urgency to transmit excess wind generation from north to south. The resulting curtailment of wind and PV generations as well as loads in the NLR and DLR simulation is shown in Table~\ref{tab:result}. It can be clearly seen that the wind curtailment is reduced significantly in the grid region of BREMEN and BERLIN with DLR. The load curtailment in the south-west region, FRANKFURT and STUTTGART is reduced by nearly half in the DLR case compared to NLR, while the load curtailment in other zones are increased because the cost function penalizes overall load curtailment.

\subsection{Line Loading and Congestion}\label{Result:II}

We select the transmission corridor BREMEN-FRANKFURT-STUTTGART to exam the real-time line loading in the scaled simulation through the year. In Fig.~\ref{fig:NLR_DLR}, the line loadings of BREMEN-FRANKFURT and FRANKFURT-STUTTGART in the NLR as well as the DLR case are shown for the simulation year 2011, the green curve represents the actual loading in the line, while the red curves are the ratings of the lines. During the winter time where wind resources are abundant, this transmission corridor is heavily loaded in the north to south direction. Clearly, DLR offers a significantly higher transmission capacity under higher winds, allowing a lot more power to be transmitted to FRANKFURT and STUTTGART. To take a closer look, the curtailment results of BREMEN and STUTTGART during December are exhibited in Fig.~\ref{fig:bremen_stuttgart}. The two set of figures illustrate that with DLR, excessive wind generations are transmitted to the south, reducing wind curtailments in BREMEN and load curtailment in STUTTGART.

\section{Discussion}


Due to limitations in the modeling method and scale, it is possible that the contribution of utilizing DLR in integrating RES is over appraised in this study. In reality, weather conditions change over spatial distance, causing the DLR to vary considerably along the transmission line. The actual transmission capacity of this line is therefore determined by the minimal DLR available over all its line sections, which is identical to the principle illustrated by Liebig's law of the minimum. An accurate estimation of a line's DLR therefore must be made based on the weather knowledge over all areas passing through by the line. The node simplification made in the benchmark zone modeling in this work ignores the unevenly distributed weather, and thus resulting in higher DLR estimations than actual values.

The rating of grid equipments must also be addressed in adopting DLR. The transmission capacity of a line is limited by grid equipments such as breakers, transformers, and converters. The power rating of these grid components must be upgraded to carry the extra power brought by DLR. Fitting the equipment rating to the merely occurred highest DLR is an costly solution due to DLR's stochastic behavior inherited from weathers. Therefore, determining a robust cap for DLR that will also account for worst-case scenarios (extreme ambient conditions in conjunction with DLR measurement errors) with a very high probability~\cite{bucher2013probabilistic}.

Due to the distribution weather conditions and power limitations from grid equipments, implementing schemes to employ DLR in power system dispatching can be challenging. It requires real-time weather monitoring over the operating area, and careful design on DLR value cap as well as equipment upgrades. The conductor temperature approach described by Foss et al. in~\cite{dlr_1990} can be a more efficient way compared to calculating DLR from weather conditions. This approach can possibly be realized by integrating temperature monitoring of transmission lines and the surrounding air mass into SCADA systems, but such developments still need considerably efforts. However, applying DLR to the most congested transmission lines is a worthwhile action, and can significantly improve the system's performance in integrating RES. As shown in Section \ref{Result:II}, employing DLR can significantly reduce load curtailments in the grid region of FRANKFURT and STUTTGART. From Fig.~\ref{fig:bench_topo}, the BREMEN-FRANKFURT-STUTTGART transmission corridor consists of two 220kV and seven 330 kV HV lines. Employ DLR on these lines can allow approximately 800 GWh more wind energy to be integrated into the grid through a whole year. 



\section{Conclusion}

The DLR model and a weather data reconstruction method have been established to calculate DLR from daily temperature records and RES feed-in series. The Germany-based economic dispatch simulation using NLR and the RES feed-in and load demand data over the year 2011 shows no renewable generation or load curtailment, showing that the German grids are capable of integrating RES generations and no need for DLR at present. However, the simulation with scaled RES generation is not merely an arbitrary hypothetical scenario but describes a probable future situation in 10 to 20 years time, especially in light of the unexpectedly rapid deployment of wind and PV units in the past. RES generation development is still growing rapidly in Germany, while all the nuclear power plants are being shut down. By that time, instead of building new transmission lines, introducing DLR might be able to solve or at least alleviate the problem.



\bibliographystyle{IEEEtran}	
\bibliography{IEEEabrv,literature}		

\end{document}